\documentclass[a4paper,12pt]{extarticle}
\usepackage{cmap}                        
\usepackage[utf8]{inputenc}             
\usepackage[english]{babel}
\usepackage[left=2cm,right=2cm,top=2cm,bottom=2cm]{geometry}  
\usepackage{amssymb, amsmath,  amsthm, amsfonts}

\theoremstyle{plain}

\newtheorem{theorem}{Theorem}

\newtheorem{corollary}{Corollary}

\theoremstyle{definition}
\newtheorem{problem}{Problem}

\begin{document}

\begin{center}
  \Large{\textbf{On supersoluble like groups of automorphisms of compact Riemann surfaces}}\footnote{The work was carried out with the financial support of the programm ``Convergence -- 2030'' (Task 1.01.4, №20262141).}\normalsize
  
  \hspace{5mm}
  
  Viachaslau I. Murashka
  
  \hspace{5mm}
  
  \small Faculty of Mathematics and Technologies of Programming,\\ Research Laboratory ``Mathematics of Hybrid Intelligence Systems'',\\ Francisk Skorina Gomel State University, Gomel,  Belarus
  
  \small mvimath@yandex.by

\medskip

\emph{Dedicated to the Professor Alexander Vasil'ev on the occasion of his 70th birthday}.

\end{center}

\begin{abstract}
  Let $G$ be a not necessary full  automorphism group of a compact Riemann surface of genus $g$. We found sharp bounds for $|G|$ in the terms of $g$ when $G$ belongs to one of the following classes of finite groups sharing similar properties with supersoluble groups: 2-nilpotent, with the Sylow tower of supersoluble type, with $\mathbb{P}$-subnormal cyclic primary subgroups, widely supersoluble, with submodular Sylow subgroups, with submodular cyclic primary subgroups and strongly supersoluble.


\end{abstract}

\section*{Introduction and the main results}

It is known that any finite group is isomorphic to a subgroup  of the full automorphism group of a compact Riemann surface. According to Hurwitz \cite{Hurwitz1893} 
a compact Riemann surface of genus $g \geq 2$ has at most $84(g - 1)$ automorphisms. Nevertheless a soluble  group of automorphisms of a compact Riemann surface of genus $g \geq 2$ has at most $48(g - 1)$ automorphisms \cite{Chetiya1981}. 

\begin{problem}[\cite{Gromadzki1988, Gromadzki1989}]\label{p1}
  Let for an integer $g\geq 2$ and a class $\mathfrak{F}$ of
finite groups, $N(g, \mathfrak{F})$ denote the order of the largest group in $\mathfrak{F}$ that a compact
Riemann surface of genus $g$ admits as a (not necessary full) group of automorphisms. Having a class $\mathfrak{F}$
of finite groups it is natural to ask for a bound for $N(g, \mathfrak{F})$ as well as for the
characterizations of those $g$ for which this bound is attained.
\end{problem}

The upper bound for $N(g, \mathfrak{F})$ as well as infinite series of $g$ for which this bound is attained  were found for the classes of all $p$-groups \cite{Zomorrodian1987},  nilpotent groups $\mathfrak{N}$ \cite{Zomorrodian1985},   supersoluble groups $\mathfrak{U}$ \cite{Gromadzki1989, Zomorrodian1990},  abelian groups $\mathfrak{A}$ \cite{Gromadzki1988}, metabelian groups \cite{Chetiya1986} and soluble groups $\mathfrak{S}$ \cite{Chetiya1981, Gromadzki1991}.  Most of the above mentioned results can be found in the survey \cite{Gromadzki1988}. Some new progress in this direction was made in \cite{Schweizer2016, Schweizer2017}.
Nevertheless, the above mentioned problem is far from its complete solution.

The paper \cite{Vasilev2010} of A.F. Vasil’ev,  T.I. Vasil’eva and V.N. Tyutyanov established  the foundations for the study of classes of finite groups of supersolvable type and their applications. In it the class $w\mathfrak{U}$ of all $w$-supersoluble finite groups was introduced. Recall that a finite group  is called $w$-supersoluble provided that every its Sylow subgroup is $\mathbb{P}$-subnormal. For applications of this class, for example see \cite{BallesterBolinches2019, Lucchini2021, Monakhov2013, Murashka2024}. In \cite{Monakhov2013} the class of finite groups with $\mathbb{P}$-subnormal cyclic primary subgroups was studied. Nowadays such groups are called $v$-supersoluble  and their class is denoted by  $v\mathfrak{U}$. The class of finite groups with submodular Sylow subgroups were studied in \cite{Monakhov2023, Vasilyev2015, Zimmermann1989}. We use $sw\mathfrak{U}$ to denote this class of groups. In \cite{Vasilyev2015} the class of all strongly supersoluble finite groups, i.e. supersoluble finite $sw\mathfrak{U}$-groups, was studied. We use $s\mathfrak{U}$ to denote the class of all such groups.  The classes $w\mathfrak{U}$, $v\mathfrak{U}$, $sw\mathfrak{U}$ and $s\mathfrak{U}$ share a lot of common properties with the class $\mathfrak{U}$ of all supersoluble groups. 
These classes are hereditary saturated formations of finite groups and are contained in the class $\mathfrak{D}$ of all finite groups with the Sylow tower of supersoluble type. Note that every $\mathfrak{D}$-group has a normal Hall $2^{\prime}$-subgroup, i.e. is $2$-nilpotent.
We use $\mathcal{N}_2$  to denote the class of all 2-nilpotent groups. It is known that above mentioned classes satisfy
$$ s\mathfrak{U}\subset\mathfrak{U}\subset w\mathfrak{U} \subset v\mathfrak{U} \subset\mathfrak{D}\subset \mathcal{N}_2 \text{ and } s\mathfrak{U}\subset sw\mathfrak{U}\subset w\mathfrak{U}. $$

In the survey \cite[Problem 3.3]{Vasilev2026} it was asked to solve Problem  \ref{p1}  for several important classes of finite group including $w\mathfrak{U}$ and $v\mathfrak{U}$. Moreover it was asked \cite[Problem 3.4]{Vasilev2026} given $g\geq 2$ and a class of groups $\mathfrak{F}$ if there exists a compact Riemann surface of genus $g$ such that its full group of automorphisms belongs to $\mathfrak{F}$. The aim of this paper is to solve three above mentioned problems for supersoluble like classes of finite groups. In our solution we will use the following principle: \emph{if $\mathfrak{F}_1\subseteq\mathfrak{F}_2$ are classes of finite groups such that $N(g, \mathfrak{F}_2)\leq h(g-1)$ and $N(g, \mathfrak{F}_1)= h(g-1)$ for infinitely many $g$, then for any class of finite groups $\mathfrak{F}$ with $\mathfrak{F}_1\subseteq\mathfrak{F}\subseteq\mathfrak{F}_2$  holds $N(g, \mathfrak{F})\leq h(g-1)$ and this bound is attained for infinitely many $g$.}
Our main result is

\begin{theorem}\label{thm1}
Let   $\mathfrak{X}$ be a class of groups and $G\in\mathfrak{X}$ be a not necessary full group of automorphisms of a compact Riemann surface $X$   of genus $g$.

If $g=2$  and $s\mathfrak{U}\subseteq\mathfrak{X}\subseteq \mathcal{N}_2$, then  $|G|\leq 24$. 
For $g\geq 3$ the following holds:

$(1)$ If $\mathfrak{D}\subseteq\mathfrak{X}\subseteq \mathcal{N}_2$, then $|G|\leq 30(g-1)$. 

$(2)$ If   $sw\mathfrak{U}\subseteq \mathfrak{X}\subseteq v\mathfrak{U}$, then $|G|\leq 21(g-1)$. 

$(3)$ Let $\mathfrak{NA}$  denote the class of all  groups with the nilpotent derived subgroup. If  $s\mathfrak{U}\subseteq \mathfrak{X}\subseteq \mathfrak{NA}\cap\mathfrak{D}$, then $|G|\leq 18(g-1)$. 

In every case $(1)$--$(3)$ the equality holds for infinitely many $g$. Moreover  if the equality holds, then $G$ is the full group of automorphisms of $X$. 
\end{theorem}

Note that $\mathcal{N}_2$ is the class of all groups whose  elements of odd order form a subgroup.

\begin{corollary}[{\cite[Theorem 7.4]{Schweizer2017}}]
Let $X$ be a compact Riemann surface of genus $g\geq 2$ and $G$ be a group of automorphisms of $X$ such that the elements of odd order in $G$ form a subgroup. Then $|G| \leq
30(g-1)$ and there are infinitely many $g$  for which
this bound is sharp.
\end{corollary}

From \cite[Definition 2, Theorem 2.10]{Vasilev2010} and \cite[Theorem C]{Vasilyev2015} it follows that $sw\mathfrak{U}\subset w\mathfrak{U}\subset v\mathfrak{U}$.

\begin{corollary}
Let $X$ be a compact Riemann surface of genus $g\geq 2$ and $G$ be a $w\mathfrak{U}$-group of automorphisms of $X$. Then $|G|\leq  21(g-1)$ for $g\geq 3$, there are infinitely many $g$  for which
this bound is sharp and $|G|\leq  24$ for $g= 2$.
\end{corollary}

In \cite{Monakhov2023} the class $v\mathfrak{U}_1$ of all finite groups with submodular cyclic primary subgroups was introduced. From \cite[Lemma 3.2, Corollary 3.5]{Monakhov2023} it follows that $sw\mathfrak{U}\subset v\mathfrak{U}_1\subset v\mathfrak{U}$.

\begin{corollary}
Let $X$ be a compact Riemann surface of genus $g\geq 2$ and $G$ be a $v\mathfrak{U}_1$-group of automorphisms of $X$. Then $|G|\leq  21(g-1)$ for $g\geq 3$, there are infinitely many $g$  for which
this bound is sharp and $|G|\leq  24$ for $g= 2$.
\end{corollary}

It is well known that $\mathfrak{U}\subset \mathfrak{NA}\cap\mathfrak{D}$. Note that $s\mathfrak{U}\subset \mathfrak{U}$ by  \cite[Theorem A]{Vasilyev2015}.

\begin{corollary}[{\cite[Theorem A]{Gromadzki1988}, \cite[Theorem 3.3]{Zomorrodian1990}}]
  A supersoluble group of automorphisms of a compact Riemann surface
of genus $g\geq 3$ has no more than $18(g - 1)$ elements. A surface of genus $g = 2$ can admit a
supersoluble group of automorphisms of order 24.
\end{corollary}

\section{Preliminaries}

\subsection{Groups of automorphisms of compact Riemann surfaces}

Recall the following well-known concepts and results on the groups of automorphisms of compact Riemann surfaces. By a group of automorphisms of a compact Riemann surface we mean any subgroup of its full group of automorphisms.
A 
Fuchsian group is a
discrete subgroup $\Gamma$ of orientation preserving isometries of the upper half plane $D$ with
hyperbolic structure. In the case when $D/\Gamma$ is a compact surface it has the presentation of the form
$$\langle x_1, \dots, x_r, a_1, b_1, \dots, a_g, b_g \mid x_1^{m_1}=\dots=x_r^{m_r}=\prod_{i=1}^{r} x_i \prod_{i=1}^{g} a_i^{-1} b_i^{-1} a_i b_i=1 \rangle.$$
A Fuchsian group with such presentation is said to have a signature $(g; m_1, \dots, m_r)$ where $g$ is called the genus and $m_i$ are called periods. The hyperbolic measure of the  associated  fundamental region for $\Gamma$ is given by
$$\mu(\Gamma) = 2\pi \left(2g - 2 + \sum_{i=1}^{r} \left(1 - \frac{1}{m_i}\right)\right).$$
For any subgroup $\Gamma_1$ of a finite index of $\Gamma$ the following Riemann -- Hurwitz formula holds
$$\mu(\Gamma_1)/\mu(\Gamma)=|\Gamma:\Gamma_1|.$$
According to the uniformization theorem, any compact Riemann surface $X$ of genus $g \geq 2$ can be expressed as the quotient $D/\Gamma$, where $\Gamma$ denotes a Fuchsian group with signature $(g; -)$, referred to as the surface group of genus $g$. 
A finite group $G$ is a group of automorphisms of $X$ if and only if there exists a Fuchsian group $\Lambda$ and a homomorphism $\theta$ from $\Lambda$ onto $G$ whose kernel is $\Gamma$. Such a homomorphism is called a surface-kernel homomorphism, and a group $G$ represented in this form is called a surface-kernel factor group of $\Lambda$. It is well known that $\theta$ is a surface-kernel homomorphism if and only if it preserves the periods of $\Lambda$.

The following result of   Macbeath \cite{Macbeath1961} plays an important role in this paper. 

\begin{theorem}[\cite{Macbeath1961}]\label{thm}
   If there is a compact Riemann surface $X$ of genus $g \geq 2$ with a group $G$ of automorphisms  with $|G|=h(g-1)$, where $h$ is a rational number with denominator dividing $g-1$, then for any positive infinite $m$ there is a compact Riemann surface of genus $g_m=m^{2g}(g-1)+1$  with a group $G_m$ of automorphisms  with $|G_m|=h(g_m - 1)$. Moreover $G$ is a factor group of $G_m$.
\end{theorem}

\subsection{Classes of finite groups}

Recall that a group is called $p$-nilpotent if it has a normal Hall $p^{\prime}$-subgroup; a class of groups is a collection of groups closed under taking isomorphic images; a class of groups is called a formation if it is closed under taking epimorphic images and finite subdirect products.
We use $Z_n$ to denote the cyclic group of order $n$. The point $(1)$ of the next theorem is the well known characterization of supersoluble groups.

\begin{theorem}\label{thm3}
  Let 
  $$1=G_0\trianglelefteq G_1\trianglelefteq\dots\trianglelefteq G_n=G$$
  be a chief series of a finite group $G$. Then
  
  $(1)$ $G$ is supersoluble iff $G/C_G(H/K)$ is an abelian group of exponent dividing $p-1$ for any chief factor $H/K$ of $G$ and $p\in\pi(H/K)$.
  
  $(2)$ \emph{\cite[Theorem 2]{Murashka2014}} $G\in v\mathfrak{U}$   iff $G/C_G(H/K)$ has   exponent dividing $p-1$ for any chief factor $H/K$ of $G$ and $p\in\pi(H/K)$.
  
  $(3)$  \emph{\cite[Theorem C]{Vasilyev2015}} $G$ has submodular Sylow subgroups iff $G/C_G(H/K)$ has  abelian Sylow subgroups of square free exponent dividing $p-1$ for any chief factor $H/K$ of $G$ and $p\in\pi(H/K)$.
  
  $(4)$  \emph{\cite[Theorem A]{Vasilyev2015}} $G$ is strongly supersoluble iff $G/C_G(H/K)$ is an abelian group of square free exponent dividing $p-1$ for any chief factor $H/K$ of $G$ and $p\in\pi(H/K)$.
  
\end{theorem}

\section{Proves of the main results}

\subsection{The proof of Theorem \ref{thm1}}

At first assume that $g\geq 3$.  

Let prove $(1)$. Note that the signature of a Fuchsian group $\Gamma$ with $\mu(\Gamma)< 2\pi/15$ can only be
$$(0;2,3,7), (0;2,3,8), (0;2,3,9), (0;2,4,5).$$
Note that the orders of the corresponding surface kernel factor groups are respectively
$$84(g-1), 48(g-1), 36(g-1), 40(g-1).$$
Hence the orders of such groups are divisible by 4. From the other hand  if a Fuchsian group $\Gamma$ of any of the previously mentioned signatures has a finite 2-nilpotent surface kernel factor group $G$, then $G$ has  the normal Hall $2^{\prime}$-subgroup    $H$  and $|G/H|=2$.
Hence $|G|$ is not divisible by $4$, a contradiction.   

Note that a Fuchsian group $\Gamma$ with $\mu(\Gamma)= 2\pi/15$ has signature $(0;2,3,10)$.  
Let $G=\langle a, b\rangle$ where
\begin{align*}
  a=&(6,15)(7,11)(8,12)(9,13)(10,14),\\
  b=&(1,6,11)(2,7,12)(3,8,13)(4,9,14)(5,10,15),\\
  (ab)^{-1}= &(1,7,2,8,3,9,4,10,5,6)(11,15,14,13,12).
\end{align*}
Then $G$ is the surface kernel factor group of $\Gamma$ of order 150. It is straightforward to check that $G\in \mathfrak{D}$.  So $N(\mathfrak{D}, 6)=30(6-1)$. Hence by Theorem \ref{thm} there is a compact Riemann surface of genus $g_{5^t}=5^{12t+1}+1$ such that it has a group $G_1$ of automorphisms with $|G_1|=30(g_{5^t}-1)$ and $G_1/H\simeq G$ for some $H\trianglelefteq G_1$. Note that $H$ is a $5$-group. So $\mathrm{O}_5(G_1/H)=\mathrm{O}_5(G_1)/H$. Since $\mathrm{O}_5(G)$ is the normal Sylow 5-subgroup of $G$, we see that $\mathrm{O}_5(G_1)$ is the normal Sylow 5-subgroup of $G_1$. Now $G_1/\mathrm{O}_5(G_1)$ is a $\{2, 3\}$-group with a Sylow tower of supersoluble type. Thus $G_1\in\mathfrak{D}$.

It means that if $\mathfrak{D}\subseteq\mathfrak{X}\subseteq \mathcal{N}_2$ and  $G\in\mathfrak{X}$ is a group of automorphisms of a compact Riemann surface $X$   of genus $g\geq 3$, then $|G|\leq 30(g-1)$. 
Assume that $|G|= 30(g-1)$. If $G$ is not a full group of automorphisms of $X$, then $G$ is a subgroup of group $H$ with $|H|\in\{84(g-1), 48(g-1), 36(g-1), 40(g-1)\}$, a contradiction with the Lagrange's theorem. Thus $(1)$ is proved.

Let prove $(2)$. Assume now that there is a surface kernel factor group $G\in v\mathfrak{U}$ of the Fuchsian group $\Gamma$  of signature $(0;2,3,10)$. Moreover assume that $G$ is a minimal order group with this property. Then $G$ is generated by $x, y$ with $o(x)=2, o(y)=3, o(xy)=10$. If $G$ is not a $\{2, 3, 5\}$-group, then $G\in v\mathfrak{U}\subset\mathfrak{D}$ has a normal Hall $\{2, 3, 5\}^{\prime}$-subgroup $H$. Hence $\Gamma$ admits $G/H$ as a surface kernel factor group, a contradiction with $|G/H|<|G|$. So $G$ is a $\{2, 3, 5\}$-group. Let $P$ be a Sylow $5$-subgroup of $G$ and $$1=P_0\leq P_1\leq\dots\leq P_n=P$$ be a part of a chief series of $G$ below it. Since $G\in v\mathfrak{U}$, we see that $G/C_G(P_i/P_{i-1})$ has exponent dividing $4$ for any $i\in\{1, \dots, n\}$ by $(2)$ of Theorem \ref{thm3}. Since the class of all groups of exponent dividing $4$ is a formation, we see that $G/\cap_{i=1}^n C_G(P_i/P_{i-1})$ has exponent dividing $4$. Now $G/C_G(P)$ is the extension of a $5$-group by a group of exponent dividing $4$ by \cite[A, Corollary 12.4(a)]{Doerk1992}. Therefore $C_G(P)$ contains every 3-element of $G$. Let $H$ be a Hall $\{3, 5\}$-subgroup of $G$ and $T$ be its Sylow 3-subgroup. Then $H\trianglelefteq G$, $P\leq H$ and $T$ is a Sylow 3-subgroup of $G$. Now $P\leq C_H(T)$, i.e. $T\trianglelefteq H$ and hence $T\trianglelefteq G$. It means that the product of 3-elements is a 3-element.
Now $xyxy=y^xy$ is a $3$-element of order 5, a contradiction. Therefore  $\Gamma$ does not admits a surface kernel factor group from $v\mathfrak{U}$.

Note that the signature of a Fuchsian group $\Gamma$ with $2\pi/15<\mu(\Gamma)< 4\pi/21$ can only be
$$(0; 2, 3, 11), (0; 2, 3, 12), (0; 2, 3, 13), (0; 2, 4, 6).$$
The groups with signatures $(0; 2, 3, 11)$ and $(0; 2, 3, 13)$ are perfect and therefore does not have soluble surface kernel factor groups.

Assume now that there is a surface kernel factor group $G\in v\mathfrak{U}$ of the Fuchsian group $\Gamma$  of signature $(0;2,3,12)$. Moreover assume that $G$ is a minimal order group with this property. Then $G$ is generated by $x, y$ with $o(x)=2, o(y)=3, o(xy)=12$. If $G$ is not a $\{2, 3\}$-group, then $G\in v\mathfrak{U}\subset\mathfrak{D}$ has a normal Hall $\{2, 3\}^{\prime}$-subgroup $H$. Hence $\Gamma$ admits $G/H$ as a surface kernel factor group, a contradiction with $|G/H|<|G|$. So $G$ is a $\{2, 3\}$-group. From $G\in\mathfrak{D}$ it follows that
the product of 3-elements is a 3-element.  Now $xyxy=y^xy$ is a $3$-element of order 6, a contradiction. Therefore  $\Gamma$ does not admits a surface kernel factor group from $v\mathfrak{U}$.

Assume now that there is a surface kernel factor group $G\in v\mathfrak{U}$ of the Fuchsian group $\Gamma$  of signature $(0;2,4,6)$. By analogy with the previous paragraph we may assume that $G$ is a minimal order group with this property and show that $G$ is a $\{2, 3\}$-group. Let $P$ be a Sylow 3-subgroup of $G$. So $P\trianglelefteq G$. Then $o(xP)=2, o(yP)=4, o(xyP)=2$. Hence $xyP=y^3xP$. It means that $G/P\simeq D_8$ is the dihedral group of order 8. So the exponent of $G/P$ is 4. Let $$1=P_0\leq P_1\leq\dots\leq P_n=P$$ be a part of a chief series of $G$ below $P$. Since $G\in v\mathfrak{U}$, we see that $G/C_G(P_i/P_{i-1})$ has exponent dividing $2$ for any $i\in\{1, \dots, n\}$ by $(2)$ of Theorem \ref{thm3}.
It means that $G/C_G(P_i/P_{i-1})$ is isomorphic to a proper quotient group of $G/P$ and hence is abelian. Therefore  $G/C_G(P_i/P_{i-1})$ is an abelian group of exponent dividing $2$ for any $i\in\{1, \dots, n\}$. It means that  $P_i/P_{i-1}$ is cyclic for any $i\in\{1, \dots, n\}$. Therefore $G$ has a chief series with  cyclic chief factors. Thus $G$ is supersoluble. From \cite{Gromadzki1989} it follows that $\Gamma$  can admit a supersolube surface kernel factor group only for $g=2$, a contradiction with $g\geq 3$.

Note that a Fuchsian group $\Gamma$ with $\mu(\Gamma)= 4\pi/21$ has signature $(0;2,3,14)$.  
Let $G=\langle a, b\rangle$ where
\begin{align*}
  a=&(8,18)(9,19)(10,20)(11,21)(12,15)(13,16)(14,17),\\
  b=&(1,8,15)(2,9,16)(3,10,17)(4,11,18)(5,12,19)(6,13,20)(7,14,21),\\
  (ab)^{-1}=&(1,15,5,19,2,16,6,20,3,17,7,21,4,18)(8,11,14,10,13,9,12).
\end{align*}
Direct calculations shows that it has a normal subgroup of order $7^2$ with quotient isomorphic to $S_3$. Hence $G/C_G(H/K)$ has abelian Sylow subgroups of square free exponent dividing $p-1$ for every $p$-chief factor $H/K$ of $G$ and $p\in\pi(G)$. Therefore $G\in sw\mathfrak{U}$ by $(3)$ of Theorem \ref{thm3}.
So $N(sw\mathfrak{U}, 15)=21(15-1)$. Hence by Theorem \ref{thm} there is a compact Riemann surface of genus $g_{7^t}=7^{30t}\cdot 14+1$ such that it has a group $G_1$ of automorphisms with $|G_1|=21(g_{7^t}-1)$ and $G_1/H\simeq G$ for some $H\trianglelefteq G_1$. 
Note that $H$ is a 7-group. Since $G/C_G(T/K)$ does not contain any proper normal 7-subgroup for any chief factor $T/K$ of $G$ below $H$ by 
\cite[A, Lemma 13.6]{Doerk1992}, $G/C_G(T/K)$ has abelian Sylow subgroups of exponent dividing $6$. Therefore $G_1\in sw\mathfrak{U}$ by $(3)$ of Theorem \ref{thm3}. 

It means that if $sw\mathfrak{U}\subseteq\mathfrak{X}\subseteq v\mathfrak{U}$ and  $G\in\mathfrak{X}$ is a group of automorphisms of a compact Riemann surface $X$   of genus $g\geq 3$, then $|G|\leq 21(g-1)$. 
Assume that $|G|= 21(g-1)$. If $G$ is not a full group of automorphisms of $X$, then $G$ is a subgroup of a group $H$ with $$|H|\in\{84(g-1), 48(g-1), 36(g-1), 40(g-1), 30(g-1), 24(g-1), 132(g-1)/5, 156(g-1)/7\}.$$ From  the Lagrange's theorem it follows that $|H|=84(g-1)$. So $H$ is a surface kernel factor group of the Fuchsian group $\Gamma$ of signature $(0; 2,3,7)$, i.e. a Hurwitz group. Note that $\Gamma=\Gamma^{\prime}$. Hence $H^{\prime}=H$. Thus $H$ has no non-trivial soluble factor groups. Since $|H:G|=4$ the action of $H$ on the rights cosets of $G$ induces homomorphism from $H$ to $S_4$ with kernel $G_H$. Since $S_4$ is soluble, we have a contradiction.    Thus $(2)$ is proved. 
 
Let prove $(3)$.
Assume now that there is the Fuchsian group $\Gamma$  with a surface kernel factor group $G\in \mathfrak{D}\cap\mathfrak{NA}$. Then 
$2\pi/15\leq\mu(\Gamma)$ by $(1)$. Let prove that $2\pi/9\leq\mu(\Gamma)$. Assume  that $2\pi/15\leq\mu(\Gamma)<2\pi/9$.
Then the signature of $\Gamma$ belongs to 
$$\{(0; 2, 3, k)\mid 10\leq k\leq 17\}\cup\{(0; 2, 4, 6), (0; 2, 4, 7), (0; 2, 5, 5), (0; 3, 3, 4)\}
$$ 
Groups with signatures $(0; 2, 3, k)$ where $ k\in\{11, 13, 17\}$ are perfect and so has no surface kernel soluble factor groups.
Using the Sylow tower property and argument similar to the proof of the case (0; 2, 3, 12), i.e. that the product of odd order elements has an odd order,  we can ruled out the signatures $(0; 2, 3, 12), (0; 2, 4, 7), (0; 2, 5, 5), (0; 3, 3, 4)$. From the 
\cite{Gromadzki1989} it follows that Fuchsian groups with signatures $(0; 2, 3, 10), (0; 2, 3, 14), (0; 2, 3, 15), (0; 2, 3, 16)$ have no surface kernel factor group with the nilpotent derived subgroup and if a Fuchsian group of signature $(0; 2, 4, 6)$  has a surface kernel factor group with the nilpotent derived subgroup, then $g=2<3$. We get a contradiction with our assumption.

Thus $2\pi/9\leq\mu(\Gamma)$. Assume that $2\pi/9=\mu(\Gamma)$. Then the signature of $\Gamma$ is $(0; 2, 3, 18)$. 
In \cite{Gromadzki1989}  it was proved that for any $g=3^n+1$, $n\geq 2$, there exists a supersoluble surface kernel factor group $G$ of $\Gamma$.
Note that $G$ is a $\{2, 3\}$-group. Hence $G$ has a normal Sylow $3$-subgroup $P$ and $G/P$ is a group of order 2.
Since $G/C_G(T/K)$ does not contain any proper normal 3-subgroup for any chief factor $T/K$ of $G$ below $P$ by 
\cite[A, Lemma 13.6]{Doerk1992}, we see that $G/C_G(T/K)$ is a group of order 1 or 2, i.e. an abelian group of a square free order and  of exponent dividing $3-1$. From $(4)$ of Theorem \ref{thm3} it follows that $G\in s\mathfrak{U}$.

It means that if $s\mathfrak{U}\subseteq\mathfrak{X}\subseteq \mathfrak{D}\cap\mathfrak{NA}$ and  $G\in\mathfrak{X}$ is a group of automorphisms of a compact Riemann surface $X$   of genus $g\geq 3$, then $|G|\leq 18(g-1)$. 
Assume that $|G|= 18(g-1)$. If $G$ is not a full group of automorphisms of $X$, then $G$ is a subgroup of group $H$ with $|H|=18k(g-1)$ for some integer $k>1$ by the Lagrange theorem, i.e. $|H|\geq 36(g-1)$. Thus the only possible case  $|H|= 36(g-1)$. It means that $H$ is a surface kernel factor group of 
the Fuchsian group of signature $(0; 2, 3, 9)$. From $|G:H|=2$ it follows that $G\trianglelefteq H$. Now every normal Hall subroup of $G$ is normal in $H$. It means that $H\in\mathfrak{D}$, a contradiction with $(1)$. Thus $(3)$ is proved.

Suppose that $g=2$.
Assume that a Fuchsian group $\Gamma$ with $\mu(\Gamma)< \pi/6$ has a 2-nilpotent surface kernel factor group $G$.
Note that the signature of $\Gamma$ can only be
$$
(0; 2, 3, 7), (0; 2, 3, 8), (0; 2, 3, 9), (0; 2, 3, 10), (0; 2, 3, 11), (0; 2, 4, 5).$$
From the   Feit-Thompson's theorem it follows that every 2-nilpotent group is soluble. Hence $G$ can not be a surface kernel factor group of $\Gamma$ with signature from $\{(0; 2, 3, 7), (0; 2, 3, 11)\}$. Any signature from $\{(0;2,3,8), (0;2,3,9), (0;2,4,5)\}$ can be ruled out the same way as in the proof of Theorem \ref{thm1}. Assume that $\Gamma$ has signature $(0; 2, 3, 10)$. Then $|G|=30$ and $G$ has a normal Hall subgroup $H$ of order 15. From the Sylow's theorem it follows that $H$ has a normal Sylow 3-subgroup. Hence $G$ also has a normal Sylow 3-subgroup. So the product of every 2-element and 3-element is a $\{2, 3\}$-element, a contradiction.  
It means that $\mu(\Gamma)\geq \pi/6$.

In \cite{Gromadzki1989} it was proved that a  Fuchsian group $\Gamma$ with signature $(0; 2, 4, 6)$ has a supersoluble surface kernel factor group $G$ of order 24. Hence it is $2$-nilpotent. From the other hand $\Gamma/\Gamma^{\prime}\simeq Z_2\times Z_2$. So $G/C_G(H/K)$ is an abelian group of exponent dividing 2 for every chief 3-factor $H/K$ of $G$ and every 2-chief factor of $G$ is central. Hence $G\in s\mathfrak{U}$ by $(4)$ of Theorem \ref{thm3}.  
It means that if $s\mathfrak{U}\subseteq\mathfrak{X}\subseteq \mathcal{N}_2$ and  $G\in\mathfrak{X}$ is a group of automorphisms of a compact Riemann surface $X$   of genus $2$, then $|G|\leq 24$.

{\small
\bibliographystyle{plain}
\bibliography{Hur}}

\begin{thebibliography}{10}

\bibitem{BallesterBolinches2019}
A.~Ballester-Bolinches, W.~M. Fakieh, and M.~C. Pedraza-Aguilera.
\newblock On products of generalised supersoluble finite groups.
\newblock {\em Mediterra. J. Math.}, 16(2):46, 2019.

\bibitem{Chetiya1981}
B.~P. Chetiya.
\newblock {On genuses of compact Riemann surfaces admitting solvable
  automorphism groups}.
\newblock {\em Indian J. PureAppl. Math.}, 12:1312--1318, 1981.

\bibitem{Chetiya1986}
B.~P. Chetiya and K.~Patra.
\newblock {On Metabelian Groups of Automorphisms of Compact Riemann Surfaces}.
\newblock {\em J. London Math. Soc.}, s2-33(3):467--472, June 1986.

\bibitem{Doerk1992}
K.~Doerk and T.~Hawkes.
\newblock {\em Finite soluble groups}.
\newblock Walter de Gruyter, Berlin -- New York, 1992.
\newblock 891 p.

\bibitem{Gromadzki1988}
G.~Gromadzki.
\newblock {Maximal groups of automorphisms of compact Riemann surfaces in
  various classes of finite groups}.
\newblock {\em RACSAM}, 82(2):267--275, 1988.

\bibitem{Gromadzki1991}
G.~Gromadzki.
\newblock {On Soluble Groups of Automorphism of Riemann Surfaces}.
\newblock {\em Can. Math. Bull.}, 34(1):67--73, 1991.

\bibitem{Gromadzki1989}
G.~Gromadzki and C.~MacLachlan.
\newblock {Supersoluble groups of automorphisms of compact Riemann surfaces}.
\newblock {\em Glasg. Math. J.}, 31(3):321--327, 1989.

\bibitem{Hurwitz1893}
A.~Hurwitz.
\newblock {Ueber algebraische Gebilde mit eindeutigen Transformationen in
  sich}.
\newblock {\em Math. Ann.}, 41:403--442, 1893.

\bibitem{Lucchini2021}
A.~Lucchini and D.~Nemmi.
\newblock {The non-$\mathfrak{F}$ graph of a finite group}.
\newblock {\em Math. Nachr.}, 294(10):1912--1921, October 2021.

\bibitem{Macbeath1961}
A.~M. Macbeath.
\newblock {On a Theorem of Hurwitz}.
\newblock {\em Glasgow Math. J.}, 5(2):90--96, 1961.

\bibitem{Monakhov2013}
V.~S. Monakhov and V.~N. Kniahina.
\newblock Finite groups with $\mathbb{P }$-subnormal subgroups.
\newblock {\em Ric. Mat.}, 62(2):307--322, 2013.

\bibitem{Monakhov2023}
V.~S. Monakhov and I.~L. Sokhor.
\newblock Finite groups with submodular primary subgroups.
\newblock {\em Arch. Math.}, 121(1):1--10, 2023.

\bibitem{Murashka2014}
V.~I. Murashka.
\newblock Properties of the class of finite groups with p-subnormal cyclic
  primary subgroups.
\newblock {\em Doklady of the National Academy of Sciences of Belarus},
  58(1):5--8, 2014.

\bibitem{Murashka2024}
V.I. Murashka.
\newblock Formations of finite groups in polynomial time:
  $\mathfrak{F}$-residuals and $\mathfrak{F}$-subnormality.
\newblock {\em J. Symbolic Comput.}, 122. --- Article~102271., 2024.

\bibitem{Schweizer2016}
A.~Schweizer.
\newblock On the exponent of the automorphism group of a compact riemann
  surface.
\newblock {\em Archiv der Mathematik}, 107(4):329--340, 2016.

\bibitem{Schweizer2017}
A.~Schweizer.
\newblock {Several types of solvable groups as automorphism groups of compact
  Riemann surfaces}, 2017.
\newblock arXiv:1701.00325 [math.CV].

\bibitem{Vasilev2026}
A.~F. Vasil'ev, V.~I. Murashka, and V.~G. Safonov.
\newblock {Problems of the theory of classes of finite groups through the prism
  of applications. I}.
\newblock {\em Problems of Physics, Mathematics and Technics}, (67):7--16,
  2026.

\bibitem{Vasilyev2015}
V.~A. Vasilyev.
\newblock Finite groups with submodular sylow subgroups.
\newblock {\em Siberian Math. J.}, 56(6):1019--1027, 2015.

\bibitem{Vasilev2010}
A.~F. Vasil’ev, T.~I. Vasil’eva, and V.~N. Tyutyanov.
\newblock On the finite groups of supersoluble type.
\newblock {\em Siberian Math. J.}, 51(6):1004--1012, 2010.

\bibitem{Zimmermann1989}
I.~Zimmermann.
\newblock Submodular subgroups in finite groups.
\newblock {\em Math. Z.}, 202(4):545--557, 1989.

\bibitem{Zomorrodian1985}
R.~Zomorrodian.
\newblock {Nilpotent automorphism groups of Riemann surfaces}.
\newblock {\em Trans. Amer. Math. Soc.}, 288:241--255, 1985.

\bibitem{Zomorrodian1987}
R.~Zomorrodian.
\newblock {Classification of p-groups of automorphisms of Riemann surfaces and
  their lower central series}.
\newblock {\em Glasgow Math. J.}, 29(2):237--244, 1987.

\bibitem{Zomorrodian1990}
R.~Zomorrodian.
\newblock {Bounds for the order of supersoluble automorphism groups of Riemann
  surfaces}.
\newblock {\em Proc. Amer. Math. Soc.}, 108:587--600, 1990.

\end{thebibliography}
\end{document}